
\input amstex
\documentstyle{amsppt}

\magnification=1200
\parskip 6pt
\pagewidth{5.4in}
\pageheight{7.0in}

\baselineskip=12pt
\expandafter\redefine\csname logo\string@\endcsname{}
\NoBlackBoxes                
\NoRunningHeads
\redefine\no{\noindent}
\redefine\qed{$\ \blacksquare$}

\define\C{\Bbb C}

\define\Z{\Bbb Z}
\define\N{\Bbb N}

\define\al{\alpha}

\define\de{\delta}
\define\la{\lambda}

\define\Om{\Omega}

\define\sub{\subseteq}

\define\st{\ \vert\ } 
\def\pr{^\prime}
\define\sh{\sharp}
\define\lan{\langle}
\define\ran{\rangle}

\redefine\ll{\lq\lq}
\redefine\rr{\rq\rq\ }
\define\rrr{\rq\rq}   

\redefine\det{\operatorname {det}}

\define\Map{\operatorname {Map}}
\redefine\dim{\operatorname {dim}}
\define\virt.dim{\operatorname {virt.dim}}
\define\point{\operatorname {point}}

\define\Fn{Fl^{(n)}}
\define\Gn{Gr^{(n)}}

\define\mnc{M_n\C}

\define\Si{\bar\Sigma} 
\define\Sii{\bar \Sigma_{(i)}} 
 
\define\Siw{\bar \Sigma_{w}}
\define\Sia{\bar\Sigma_{\bold a}} 
 
\define\Ci{\bar C_{(i)}} 
\define\Cw{\bar C_{w}} 
\define\Ca{\bar C_{\bold a}} 

\define\x{X}
\define\y{Y}
\define\q{Q}

\topmatter
\title Quantum cohomology and the periodic Toda lattice
\endtitle
\author Martin A. Guest and Takashi Otofuji
\endauthor
\abstract
We describe a relation between the periodic one-dimensional Toda lattice 
and the quantum cohomology of the  periodic flag manifold 
(an infinite-dimensional K\"ahler manifold).  This generalizes a 
result of Givental and Kim relating the open 
Toda lattice and the quantum cohomology of the finite-dimensional flag 
manifold. We derive a simple and explicit 
\ll differential operator formula\rr for the necessary quantum products, 
which applies both to the finite-dimensional and to the infinite-dimensional
situations.
\endabstract
\endtopmatter

\document

\subheading{Introduction}

The quantum cohomology of the full flag manifold $F_n$ of $SU_n$ is
known to be related to an integrable system, the
open one-dimensional Toda lattice.  This relation was established in 
\cite{Gi-Ki}, and a rigorous framework for the calculations
was developed  in \cite{Ci1}, \cite{Ki1}, and \cite{Lu}, building on earlier
fundamental work in quantum cohomology.
We shall give --- in the spirit of \cite{Gi-Ki} --- an analogous relation
between the quantum cohomology of the periodic flag manifold
$\Fn$ and the periodic one-dimensional Toda lattice.

Such an extension to the periodic case is perhaps not unexpected, 
but we feel that it is worth noting, for two
reasons.  First, there are several new features of the quantum
cohomology of the periodic flag manifold $\Fn$, the most obvious one 
being
that $\Fn$ is an {\it infinite-dimensional} K\"ahler manifold.
Second, very few concrete examples of this phenomenon are known 
(cf. section 2.3 of \cite{Au}).  Indeed, 
the full flag manifold $F_n$ seems to be the only example
so far, together with its generalization
\footnote{
Some comments on the case
of partial flag manifolds and their relation with Toda
lattices can be found in section 5 of \cite{Gi1}.}
$G/B$ which was accomplished in \cite{Ki2}.  
Now, $\Fn$ is an infinite-dimensional flag manifold
(of the loop group $LSU_n$), and is therefore a close relative
of this family.  However, the  periodic one-dimensional Toda lattice
is more complicated than the open one; for example its
solutions generally involve theta functions, whereas those of the
open Toda lattice are rational expressions of exponential functions.

The open one-dimensional Toda lattice is a 
(nonlinear) first-order differential equation
$$
\dot L_n(t) = [L_n(t),M_n(t)]
$$ 
where $L_n$ is the tri-diagonal matrix
$$
L_n=
\eightpoint{
\pmatrix
\x_1 & \q_1     & 0       & \cdots & \cdots      & \cdots      & 0  
\\
-1       & \x_2 & \q_2     & 0      & \cdots      & \cdots      & 0     
\\
0       & -1       & \x_3 & \q_3    & 0           & \cdots      & 0     
\\
\vdots  & \ddots  & \ddots  & \ddots & \ddots      & \ddots      & \vdots
\\
0       & \cdots  & 0       & -1      & \x_{n-2} & \q_{n-2}     & 0     
\\
0       & \cdots  & \cdots  & 0      & -1           & \x_{n-1} &
\q_{n-1}\\
0       & \cdots  & \cdots  & \cdots & 0           & -1           & \x_n
\endpmatrix
}
$$
and $M_n$ is a certain modification of $L_n$.
Here, $\x_1,\dots,\x_n$ and $\q_1,\dots,\q_{n-1}$ are functions of a
real variable $t$ with $\q_i<0$, and we assume that $\x_1+\cdots+\x_n=0$.  Let
$$
\det(L_n+\mu I) =  O_n  =  \sum_{i=0}^n O^i_n \mu^i.
$$
Then the polynomials $O^0_n,O^1_n,\dots,O^{n-1}_n$ in
$\x_1,\dots,\x_n$ and $\q_1,\dots,\q_{n-1}$ are
\ll the conserved quantities\rr of the Toda lattice,
which give rise to its integrability.  (For further explanation
of Toda lattices we refer to \cite{Ol-Pe}, \cite{Pe}, \cite{Re-Se}.)  
The result of \cite{Gi-Ki}
is that the (small) quantum cohomology algebra of 
$$
F_n=
\{ E_1\sub E_2\sub \cdots \sub E_n = \C^n \st
E_i \text{ is an $i$-dimensional linear subspace of $\C^n$ }\}
$$
is
$$
QH^\ast F_n \ \cong\ 
\C[\x_1,\dots,\x_n,\q_1,\dots,\q_{n-1}]/
\lan O^0_n,O^1_n,\dots,O^{n-1}_n \ran,
$$
where $\x_1,\dots,\x_n,\q_1,\dots,\q_{n-1}$ are regarded now
as indeterminates.
In other words, the conserved quantities of the 
open one-dimensional Toda lattice are precisely the defining
relations for the quantum cohomology algebra of $F_n$.  This
remarkable fact has been explored in a number of very
interesting papers (such as \cite{Gi2}, \cite{Ki2}, \cite{Ko1},
\cite{Ko2}, \cite{Fo-Ge-Po}).

The periodic one-dimensional Toda lattice is a differential equation
of the form
$$
\dot \Cal L_n(t) = [\Cal L_n(t),\Cal M_n(t)]
$$ 
where $\Cal L_n$ is the matrix
$$
\Cal L_n=
\eightpoint{
\pmatrix
\x_1 & \q_1     & 0       & \cdots & \cdots      & \cdots      & -z    
\\
-1       & \x_2 & \q_2     & 0      & \cdots      & \cdots      & 0     
\\
0       & -1       & \x_3 & \q_3    & 0           & \cdots      & 0     
\\
\vdots  & \ddots  & \ddots  & \ddots & \ddots      & \ddots      & \vdots
\\
0       & \cdots  & 0       & -1      & \x_{n-2} & \q_{n-2}     & 0     
\\
0       & \cdots  & \cdots  & 0      & -1           & \x_{n-1} &
\q_{n-1}\\
\q_n/z       & \cdots  & \cdots  & \cdots & 0           & -1           & \x_n
\endpmatrix
}
$$
and where $z$ is a \ll spectral parameter\rr in 
$S^1=\{ z\in\C \st \vert z\vert = 1 \}$. Thus, $\Cal L_n$ may be 
interpreted as a 
function of the real variable $t$ with values in the loop algebra
$\Map(S^1,\mnc)$.
The variables  $\x_1,\dots,\x_n$ and $\q_1,\dots,\q_{n}$ here are functions of a
real variable $t$ with $\q_i<0$, and we assume that $\x_1+\cdots+\x_n=0$ and that
$\q_1\q_2\cdots \q_n$ is constant.   Let
$$
\det(\Cal L_n+\mu I) =  P_n  =  \sum_{i=0}^n P^i_n \mu^i + A_n\frac 1z + 
B_n z
$$
where $P^k_n$, $A_n$, $B_n$ are polynomials in
$\x_1,\dots,\x_n$ and $\q_1,\dots,\q_{n}$.  The
$P^0_n,P^1_n,\dots,P^{n-1}_n$  are
\ll the conserved quantities\rr of the periodic Toda lattice.

The loop group $LSU_n=\Map(S^1,SU_n)$ plays an analogous role here
to that of the group $SU_n$ for the open Toda lattice, and
the periodic flag manifold $\Fn$ is analogous to $F_n$
(it is a complete flag manifold for an affine Kac-Moody group).
For a precise definition of $\Fn$ we refer to section 8.7
of \cite{Pr-Se}; we just remark that it is related to the 
Grassmannian model $Gr^{(n)}$ of the based loop group
$\Om SU_n\cong LSU_n/SU_n$ as follows:
$$
\Fn \ =\ 
\{ W_0\sub W_1\sub \cdots \sub W_n \st
W_i \in Gr^{(n)},
\virt.dim W_i = i-n,
\la W_n = W_0 \}.
$$
Here, $\Gn$ is a certain subspace of the Grassmannian
of all linear subspaces of the Hilbert space 
$$
H = L^2(S^1,\C^n)=\bigoplus_{i\in \Z} \la^i \C^n,
$$
and $\la W_n$ denotes the result of applying the linear
\ll multiplication operator\rr $\la$ (of $H$) to $W_n$.  
The virtual dimension is defined by
$\virt.dim W = \dim (W \cap H_-) -\dim (W^\perp \cap H_+)$,
where $H_+ = \oplus_{i\ge 0} \la^i \C^n$, 
$H_- = \oplus_{i< 0} \la^i \C^n$.

Let $H^\sh\Fn$ denote the subalgebra of the cohomology algebra
$H^\ast\Fn$ which is generated by $H^2\Fn$.
Let $QH^\sh\Fn$ denote the subalgebra of the quantum cohomology algebra
$QH^\ast\Fn$ which is generated by $H^2\Fn$.    Then our result is:
$$
QH^\sh \Fn \ \cong\ 
\C[\y_1,\dots,\y_n,\q_1,\dots,\q_{n}]/
\lan P^0_n,P^1_n,\dots,P^{n-1}_n \ran,
$$
where $\y_1,\dots,\y_n$ are related to $\x_1,\dots,\x_n$ by $\x_i=\y_i-\y_{i-
1}$
and $\y_0=\y_n$
(the precise nature of $\y_1,\dots,\y_n$ will be made clear later).

We refrain from calling this a \ll theorem\rrr, as it depends on two
provisional assumptions which we shall not attempt to justify in this paper.
These are (1) that a rigorous definition of  $QH^\sh\Fn$ is
possible, and (2) that $QH^\sh\Fn$ and 
$H^\sh\Fn \otimes \C[\q_1,\dots,\q_{n}]$ are isomorphic as
$\C[\q_1,\dots,\q_{n}]$-modules.  
Regarding (1), we have little doubt that an
appropriate definition can be given, for example
as in \cite{Be}, \cite{Ci1},
using \ll quantum Schubert calculus\rrr. 
Assumption (2) may be avoided, as we shall explain at the end of the paper.
Our calculation
is quite short, and it gives simultaneously another proof of the
result of Givental and Kim for $F_n$ 
(where assumptions (1) and (2) are unnecessary).

To conclude this introduction, we comment on two special
features of $QH^\ast \Fn$  which are not present in the case of $QH^\ast 
F_n$:

\no(i) The space $\Fn$ --- and the space
of rational curves in $\Fn$ of fixed degree ---  is infinite-dimensional. 
On the other hand, the space of rational curves of fixed
degree in $\Fn$ {\it which
intersect a fixed finite-dimensional subvariety} is finite-dimensional.
(This is an observation of \cite{At}.)  It is this property which
is primarily responsible for the existence of the quantum cohomology of 
$\Fn$.
An alternative manifestation of this property is that the first Chern
class of $\Fn$ is {\it finite} (see \cite{Fr}).

\no(ii) Because of (i), Poincar\acuteaccent e duality 
is not immediately available for $\Fn$.  However, as a substitute, we use 
the existence of dual Birkhoff and Bruhat cells in $\Fn$ (see \cite{Pr-Se}).
Bruhat cells are finite-dimensional and their closures represent 
a basis for the
homology classes of $\Fn$; Birkhoff cells are finite-codimensional
and their closures represent
a basis for the  the cohomology classes.  These play the role
of Schubert varieties and \ll dual\rr Schubert varieties in $F_n$.

Acknowledgements: The first author is grateful for financial
support from the NSF (USA), the NCTS (National Tsing Hua University, 
Taiwan), and Hull University (UK). 
The authors thank Augustin-Liviu
Mare for several helpful comments.

\subheading{\S 1 The periodic flag manifold } 

We shall review some facts concerning $\Gn$ and $\Fn$ from chapter 8 of 
\cite{Pr-Se},
and establish some additional notation. Recall that 
$\Gn$ has a line bundle $\det\Cal W$, which may be considered as the \ll 
top exterior power\rr of the
tautologous bundle $\Cal W$. (The fibre of $\Cal W$ over $W\in\Gn$ is 
$W$ itself.) Similarly, $\Fn$ has
tautologous bundles $\Cal W_i$ and associated line 
bundles $\det\Cal W_i$. 

\proclaim{Definition 1.1}

\no (1) $y_i=-c_1 \det \Cal W_i \in H^2 \Fn$ 

\no (2) $x_i=-c_1 \Cal W_i / \Cal W_{i-1} = y_i-y_{i-1} \in H^2 \Fn$. 
\endproclaim

\no The bundles $\Cal W_i / \Cal W_{i-n}$ are known to be topologically 
trivial, so we have $y_0=y_n$ and
$x_0=x_n$. Since $\Cal W_n / \Cal W_{0}$ is topologically
equivalent to $\oplus_{i=1}^n \Cal W_i / \Cal W_{i-1}$, it follows
that the elementary symmetric functions of $x_1,\dots,x_n$ are zero. 

Let $\C[\y_1,\dots,\y_n]$ be the algebra of complex polynomials in
certain variables $\y_1,\dots,\y_n$.  The \ll classical evaluation map\rr
$$ 
ev_c:\C[\y_1,\dots,\y_n] \to H^\sh \Fn,\quad \y_i \mapsto y_i
$$
is by definition an epimorphism, and we shall investigate its kernel.
For this we need the elementary symmetric polynomials in 
$\x_1,\dots,\x_n$,
where $\x_i=\y_i-\y_{i-1}$ $(i \in \Z/n\Z)$:

\proclaim{Definition 1.2}
$S_n = \sum_{i=0}^n S^i_n \mu^i = (\x_1+\mu)\cdots(\x_n+\mu)$. 
\endproclaim

\no We shall show that
the $S^i_n$ are generators of the kernel of $ev_c$, i.e.
the relations defining the algebra $H^\sh \Fn$.

\proclaim{Proposition 1.3} We have
$H^\sh \Fn \cong \C[\y_1,\dots,\y_n]/
\lan S^0_n,S^1_n,\dots,S^{n-1}_n \ran$,
the isomorphism being induced by $ev_c$. \endproclaim

\demo{Proof} Since the bundle $\Cal W_n / \Cal W_{0}$ is topologically 
trivial, 
the map $\pi_n:\Fn\to \Gn$ given by $\pi(W_0\sub W_1\sub \cdots \sub 
W_n)= W_n$ 
defines a trivial bundle over the \ll identity component\rr $\Gn_{\text 
id}$ of 
$\Gn$ consisting of
subspaces of virtual dimension zero. Now, $\Gn_{\text id}$ is homotopy 
equivalent  to $\Om SU_n=L SU_n/SU_n$ (see \cite{Pr-Se}), 
and its cohomology is well known (see \cite{Bo}). 
The fibre of the bundle is $F_n$. Hence 
$H^\sh \Fn \cong H^\ast F_n \otimes H^\sh \Om SU_n$, and this leads to 
the stated 
result.
\qed
\enddemo

We shall make use of the Birkhoff \ll cells\rr $\Sigma_{\bold a}$ and
the Bruhat cells $C_{\bold a}$  of the Grassmannian $\Gn$, which were
introduced in section 8.4 of \cite{Pr-Se}.  They are indexed
by elements $\bold a$ of $\Z^n$.  The closures $\Ca$
of the Bruhat cells are finite-dimensional projective algebraic
varieties, and their fundamental homology classes form a system
of additive generators for $H_\ast \Gn$.  There is a duality between the
$\Ca$'s and the $\Sia$'s which is analogous to the duality between
\ll opposite\rr Schubert decompositions of a finite-dimensional
Grassmannian.  This may be expressed in terms of intersections --- see
Theorem 8.4.5 of \cite{Pr-Se}.  The
finite-codimensional varieties $\Sia$ can be considered as
representatives of a system of additive generators for $H^\ast \Gn$.

For example, if $\Gn_{\text id}$ denotes the component of $\Gn$ 
consisting of
subspaces of virtual dimension zero as in the proof above, 
and $\Cal W_{\text id}$ denotes the
restriction of the bundle $\Cal W$ to $\Gn_{\text id}$, then 
the cohomology class $y=-c_1 \det \Cal W_{\text id}
\in H^2 \Gn_{\text id}$ corresponds to the
unique Birkhoff variety of codimension one in $\Gn_{\text id}$,
in the sense that the latter is the zero set of a holomorphic section of
$\det \Cal W^\ast_{\text id}$ (see section 7.7 of \cite{Pr-Se}).
This variety, which we shall denote by $\Si$, is given explicitly by
$$
\Si=\{ W\in \Gn_{\text id} \st \dim W\cap H_- \ge 1 \}.
$$
The dual Bruhat variety $\bar C$ is given explicitly by
$$
\bar C=\{ W\in \Gn_{\text id} \st 
\C e_2\oplus\cdots\oplus\C e_{n} \oplus \la H_+
\sub W \sub
\C \la^{-1} e_n \oplus H_+
 \},
$$
where $e_1,\dots,e_n$ is the standard basis of $\C^n$.

Birkhoff varieties $\Siw$ and Bruhat varieties $\Cw$ for the
periodic flag manifold $\Fn$ were defined in section 8.7 of \cite{Pr-Se}
in a similar way.  They are indexed by elements $w$ of the affine
Weyl group of $L SU_n$.  This time there are $n$ Birkhoff varieties
of codimension one, corresponding to the additive generators
$y_i$ of $H^2 \Fn$, namely
$$
\Sii = \{ (W_0\sub W_1\sub \cdots \sub W_n)\in\Fn  \st
\dim W_i \cap (H_-\oplus \C e_{1}\oplus\cdots\oplus \C e_{n-i}) \ge 1 \}.
$$
The dual Bruhat varieties are
$$
\Ci
=\{ (W_0\sub W_1\sub \cdots \sub W_n)\in\Fn  \st
W_k= \C e_{n-k+1}\oplus\cdots\oplus\C e_{n} \oplus  \la H_+
\text{ for } k\ne i \}.
$$

The inclusion of $\Ci(\cong \C P^1)$ in $\Fn$ defines a rational curve 
$f_i$,
and the homotopy classes 
$$
q_i = [f_i]
$$
form an additive basis of $\pi_2 \Fn \cong H_2 \Fn \cong \Z^n$.  
As usual in the construction of quantum cohomology, we shall in future
use  multiplicative notation $q^{\al} = q_1^{\al_1}\cdots  q_n^{\al_n}$
for an element $\al=\al_1 q_1+\cdots + \al_n q_n$ of $\pi_2 \Fn$, i.e.
instead of $\al$ we use the corresponding additive generator $q^{\al}$
of the group algebra of $\pi_2 \Fn$.  

It is easy to show (by considering the bundle
$\pi_n:\Fn\to\Gn$) that a homotopy class $q^{\al}$ contains a
rational curve only if $\al\ge 0$, i.e. $\al_i\ge 0$ for
all $i$.

\subheading{\S 2 Computations of quantum products } 

Our computations of quantum products for $\Fn$ are based on the 
existence of a Gromov-Witten invariant
$$
\lan \Si_{w_1} \vert \Si_{w_2} \vert \bar C_{w_3} \ran_{q^{\al}}. $$
This may be defined --- naively --- as the number of rational curves $f$ 
in the homotopy class $q^{\al}$ such that $$
f(0)\in g_1\Si_{w_1},\quad
f(1)\in g_2\Si_{w_2},\quad
f(\infty)\in g_3\bar C_{w_3},
$$
where $g_1,g_2,g_3$ are \ll general\rr elements of the loop group 
$LSU_n$. 
As stated in the introduction, we shall {\it assume} that such invariants 
are well defined, and that they
give rise to a commutative associative \ll quantum product\rr operation 
$\circ$ on 
$H^\ast \Fn \otimes
\C[q_1,\dots,q_n]$, through the following standard procedure: 

\proclaim{Definition 2.1} For
$u=[\Si_{w_1}], v=[\Si_{w_2}]$ in $H^\ast \Fn$, let $$
u\circ v = \sum_{\al\ge 0} (u\circ v)_{\al} q^{\al}, $$
where $(u\circ v)_{\al} \in H^\ast \Fn$ is determined via its (Kronecker) 
products
by
$$
\lan (u\circ v)_{\al}, [\Cw]\ran =
\lan \Si_{w_1} \vert \Si_{w_2} \vert \bar C_{w} \ran_{q^{\al}}, $$
for all $w$ in the affine Weyl group of $LSU_n$.
We denote by $QH^\ast\Fn$ the algebra with underlying 
$\C[q_1,\dots,q_n]$-module
$H^\ast \Fn \otimes\C[q_1,\dots,q_n]$ and product operation $\circ$.
\endproclaim

\no We assume further that $\circ$ is a deformation of the 
cup product in cohomology in the sense that $(u\circ v)_{0} = uv$, 
and that $\circ$ respects the grading defined in the usual way
by
$$
\vert x q^{\al} \vert = \vert x \vert + \lan c_1 \Fn, q^{\al}\ran
$$
(where $x\in H^\ast \Fn$). It follows from \cite{Fr} that
$\lan c_1 \Fn, q^{\al}\ran = 4\sum_{i=1}^n \al_i$, and it is easy
to check that this is the dimension of the space of basepoint
preserving rational curves in the homotopy class $q^{\al}$.
We obtain $\vert (u\circ v)_{\al}\vert = 
\vert u \vert + \vert v \vert - 4\sum_{i=1}^n \al_i$. 

Note that we are assuming, in particular, that the ordinary cup product
is given by intersections of (general translates of) Bruhat
and Birkhoff varieties.  We could not find a direct statement of this
in the literature, but it appears to be known (see
\cite{Ca}, \cite{Gu}, \cite{Ha}, \cite{Ko-Ku}).

The following useful lemma says that, for a quantum product 
of the form $y_i^m \circ v$, each nonzero term $(y_i^m\circ v)_{\al} 
q^{\al}$ in the \ll quantum deformation\rr must
be divisible by $q_i$. 

\proclaim{Lemma 2.2} Let $i\in \{ 1,\dots,n\}$, $m\in \N$. 
Let $v \in H^\ast\Fn.$
Write $y_i^m\circ v = y_i^m v + \sum_{\al> 0} (y_i^m\circ v)_{\al} q^{\al}$ 
(as above). If $(y_i^m \circ
v)_{\al}\ne 0$, then $\al_i \ge 1$. \endproclaim

\demo{Proof} The cohomology class $y_i^m$ may be represented by a 
variety of the form
$$
\Sii^m=g_1\Sii \cap \cdots \cap g_m\Sii
$$
where $g_1,\dots,g_m$ are suitable elements of $LSU_n$; this is a 
subset of $\Fn$ consisting of elements $W_0\sub W_1\sub\cdots\sub 
W_n$ 
for which $W_i$ (and only $W_i$)
satisfies a certain condition. We claim that
$$
\lan \Sii^m
\st
\Si\pr
\st
\bar C\pr
\ran_{q^{\al}}\ne 0
\implies
\al_i\ge 1,
$$
for any Birkhoff variety $\Si\pr$ and any Bruhat variety $\bar C\pr$. 

If this assertion is false, there is a (nonzero) finite number of rational 
curves
$$
f = (W_0\sub W_1\sub\cdots\sub W_n):\C P^1\to \Fn $$
in the homotopy class $q^{\al}$, with $\al_i=0$, such that $$
f(0)\in \Sii^m,\quad
f(1)\in g_1\Si\pr,\quad
f(\infty)\in g_2\bar C\pr
\quad
(\text{for some } g_1,g_2\in LSU_n).
$$
Since $\al_i=0$, $W_i$ is constant. But then we obtain a continuous 
family 
of rational curves with the same properties, by pre-composing with 
fractional linear transformations $\xi$
such that
$\xi(0)=z\in \C P^1-\{1,\infty\}$,
$\xi(1)=1$,
$\xi(\infty)=\infty$.
This is a contradiction.
\qed
\enddemo

We shall be interested in the \ll quantum versions\rr $QS^i_n$ 
of the relations $S^i_n$ of $H^\sh\Fn$.  These will be the relations for
the algebra $QH^\sh\Fn$, which is defined analogously to $H^\sh\Fn$
as the subalgebra of $QH^\ast\Fn$ generated by $H^2\Fn$
(but see assumptions (1) and (2) of the introduction).
By an argument of Siebert and Tian (Theorem 2.2 of \cite{Si-Ti}),
$QS^i_n$ is obtained by suitably modifying $S^i_n$. 
To explain this modification,  the following two facts are needed:

\no{(1)} Any quantum product can be expressed as a 
linear combination of classical products with coefficients in 
$\C[q_1,\dots,q_n]$.

\no{(2)} Any classical product can be expressed as a 
linear combination of quantum products with coefficients in 
$\C[q_1,\dots,q_n]$.

\no The first is obvious from the definition of quantum product, and the
second may be proved by an induction argument, bearing in mind that
the degree of a quantum product is the sum of the degrees of the 
individual factors.

Let $\C[\y_1,\dots,\y_n,\q_1,\dots,\q_n]$ 
be the algebra of complex polynomials in
certain variables $\y_1,\dots,\y_n,\q_1,\dots,\q_n$.
We define the quantum evaluation map
$$
ev_q: \C[\y_1,\dots,\y_n,\q_1,\dots,\q_n] \to QH^\sh \Fn
$$
as the algebra epimorphism that sends $\y_i$ to  
$y_i$ and $\q_i$ to $q_i$ in
$QH^\sh \Fn$. Via the module 
identification $QH^\sh \Fn \cong H^\sh \Fn
\otimes\C[q_1,\dots,q_n]$,
$ev_q$ can be regarded as the map which evaluates all the quantum 
products
in a \ll quantum polynomial\rr involving 
$y_1,\dots,y_n,q_1,\dots,q_n$.
The classical evaluation map $ev_c$ extends to an algebra 
epimorphism
$$
ev_c:\C[\y_1,\dots,\y_n,\q_1,\dots,\q_n] \to H^\sh \Fn 
\otimes\C[q_1,\dots,q_n].
$$

In general, $ev_c$ and $ev_q$ do not coincide, of course. But
it follows from (1) and (2) above that, for any  polynomial $R$, there is
a  (not in general unique) polynomial $\Cal R$ such that 
$ev_c R = ev_q \Cal R$.  Our  main
computational result is that there is a simple algebraic formula for 
the polynomial $\Cal  R = QS^i_n$
in terms of the polynomial $R = S^i_n$. It will be convenient to express 
this in terms of the differential operators 
$$
\delta_i = Id - \q_i\frac{\partial^2}{\partial \x_i \partial \x_{i+1}},\quad
D_i = Id + \q_i\frac{\partial^2}{\partial \x_i \partial \x_{i+1}},\quad 
i\in \Z/ n\Z.
$$
These operators commute (since they have constant coefficients).

Denote by $V$ the $\C[\q_1,\dots,\q_n]$-submodule of 
$\C[\y_1,\dots,\y_n,\q_1,\dots,\q_n]$ that is 
generated by elements of the form
$\x_{i_1} \cdots \x_{i_k}$, $1\le i_1 < \dots < i_k \le n$.

\proclaim{Proposition 2.3}
On $V$ we have

\no{(1)}
$ev_q=
ev_c \delta_n \delta_{n-1}\cdots \delta_1 $. 

\no{(2)}
$ev_c=
ev_q D_n D_{n-1}\cdots D_1 $.
\endproclaim

We shall postpone the proof of Proposition 2.3 for a moment.
Part (2)  gives our explicit formula for $\Cal R$ in terms
of $R$, namely $\Cal R=D_n D_{n-1}\cdots D_1  R$. Applying this to the
relations $R=S^i_n$ we obtain the required quantum modifications:

\proclaim{Definition 2.4}
$QS_n = \sum_{i=0}^n QS^i_n \mu^i = D_n D_{n-1}\cdots D_1 S_n$. 
\endproclaim

\proclaim{Corollary 2.5}  Subject to the validity of 
assumptions (1) and (2) of
the introduction, we have
$QH^\sh \Fn \cong \C[\y_1,\dots,\y_n,\q_1,\dots,\q_n]/
\lan QS^0_n, QS^1_n,\dots,QS^{n-1}_n \ran$,
the isomorphism being induced by $ev_q$.
\endproclaim

For example, $S^{n-1}_n=\sum_{i}\x_i = QS^{n-1}_n$
and
$$
S^{n-2}_n=\sum_{i<j}\x_i \x_j, \quad
QS^{n-2}_n=\sum_{i<j}\x_i \x_j +\sum_{i} \q_i.
$$
The  relation $QS^{n-2}_n$ 
corrresponds to the quantum multiplication formula
$\sum_{i<j}x_i\circ x_j = \sum_{i<j}x_i x_j - \sum_{i} q_i$.
The formula can be established by showing that 
$x_i\circ x_{i+1} = x_i x_{i+1} -  q_i$ for all $i$ and
$x_i\circ x_{j} =x_i x_{j}$ when $j>i+1$.  This, and its
generalization to products of the form 
$x_{i_1}\circ \cdots \circ x_{i_k}$, is the basis of our proof
of Proposition 2.3.

\demo{Proof of Proposition 2.3}
We have
$
D_i \delta_i = Id
$ on $V$,
since
$D_i \delta_i =
Id-\q_i^2 \frac{\partial^4}{\partial \x_i^2 \partial \x_{i+1}^2}$, 
and the second term vanishes on $V$. 
If we assume (1), then we have 
$$
\align
ev_q(D_n D_{n-1}\cdots D_1 \x_{i_1} \cdots \x_{i_k})  
&=
ev_c(\delta_n \delta_{n-1}\cdots \delta_1 D_n D_{n-1}\cdots D_1 
\x_{i_1} \cdots
\x_{i_k}) \\ 
&= ev_c(\x_{i_1} \cdots \x_{i_k}).
\endalign
$$
So (2) is a consequence of (1).

To prove (1), it suffices to show that the quantum
product 
$$ x_{i_1}\circ \cdots \circ x_{i_k}, \quad 1\le i_1
< \dots < i_k \le n 
$$ 
is obtained by replacing (in any order) each
occurrence of $\x_i \x_{i+1}$ by $\x_i \x_{i+1} - \q_i$ and then
applying the classical evaluation map $ev_c$. 
For example, $X_1X_2X_3$ becomes 
$X_1X_2X_3-X_3Q_1 -X_1Q_2 - X_2Q_3$.  In
terms of $y_1,\dots, y_n$ (using  $x_i=y_i-y_{i-1}$) this is equivalent to: 
\enddemo

\proclaim{Lemma 2.6} Let $a$ and $b$ be nonnegative integers.
Then
$$
(y_{i_1}\circ y_{i_1})
\circ
\cdots
\circ
(y_{i_a}\circ y_{i_a})
\circ
y_{j_1} \circ \cdots \circ y_{j_b} =
(y_{i_1}^2+q_{i_1}) \cdots(y_{i_a}^2+q_{i_a}) y_{j_1} \cdots y_{j_b}
$$
provided that 
all the indices in this expression 
are distinct and no two of $i_1,\dots,i_a$ are
consecutive (mod $n$).
\endproclaim

\demo{Proof} 
We use induction on $a+b$. For $a+b\le 1$, the only nontrivial case 
to be established is
$
y_i\circ y_i = y_i^2 + q_i.
$
By Lemma 2.2, each term in the quantum deformation of $y_i\circ y_i$ 
must contain
$q_i$. Hence $y_i\circ y_i = y_i^2 + \la q_i$, where (by definition of the 
quantum
product) we have $\la = \lan \Sii \vert \Sii \vert \point \ran_{q_i}$. This 
is
evaluated by counting rational curves $f = (W_0\sub W_1\sub\cdots\sub 
W_n)$ with
$W_k(z)$ constant if $k\ne i$. But this is essentially the Gromov-Witten 
invariant
$\lan \point
\vert \point \vert \point \ran_1$ in the quantum cohomology of $\C P^1$, 
so 
$\la =1$. 

Now we proceed to the inductive step.  By the previous paragraph, we have 
$$
(y_{i_1}\circ y_{i_1})
\circ
\cdots
\circ
(y_{i_a}\circ y_{i_a})
\circ
y_{j_1} \circ \cdots \circ y_{j_b} =
(y_{i_1}^2 + q_{i_1})
\circ
\cdots
\circ
(y_{i_a}^2 + q_{i_a})
\circ
y_{j_1} \circ \cdots \circ y_{j_b}
$$
It suffices to show that
$$
y_{i_1}^2\circ\cdots \circ y_{i_a}^2\circ y_{j_1} \circ \cdots \circ 
y_{j_b}
=
y_{i_1}^2 \cdots y_{i_a}^2 y_{j_1} \cdots y_{j_b},
$$
i.e. that products of the form $y_{i_1}^2\circ\cdots \circ y_{i_a}^2\circ 
y_{j_1} \circ \cdots \circ
y_{j_b}$ have no \ll quantum deformation\rrr.
(If $a=0$ this is the same as the statement that we wish to
prove; if $a>0$ it implies the required statement, because
of the induction hypothesis.)

We shall give a separate inductive argument for the last statement.
For $a+b\le 1$ there is nothing to prove. For the inductive step,
we consider first the case where $b>0$. 
By the inductive hypothesis, we have 
$$
\align
y_{i_1}^2\circ\cdots \circ y_{i_a}^2\circ y_{j_1} \circ \cdots \circ 
y_{j_b}
&=
y_{i_p}^2\circ
(y_{i_1}^2 \cdots \hat y_{i_p}^2 \cdots y_{i_a}^2 y_{j_1} \cdots 
y_{j_b})\\ &=
y_{j_q}\circ
(y_{i_1}^2 \cdots y_{i_a}^2 y_{j_1} \cdots \hat y_{j_q} \cdots y_{j_b}). 
\endalign
$$
Applying Lemma 2.2, we see that each term of the quantum 
deformation of the left hand side  must contain $q_{i_1}\cdots
q_{i_a} q_{j_1} \cdots q_{j_b}$. The former has degree $4a+2b$, and the 
latter has
degree $4a+4b$. Since $b>0$, this means that there is in fact no quantum
deformation.

It remains to prove the inductive step in the case where $b=0$.  By the
inductive hypothesis and Lemma 2.2 again, we have
$$
\align
y_{i_1}^2\circ\cdots \circ y_{i_a}^2 &=
y^2_{i_r} \circ (y_{i_1}^2\cdots \hat y_{i_r}^2 \cdots y_{i_a}^2)
\ \text{for any $r$}\\
&=y_{i_1}^2\cdots  y_{i_a}^2 + \la q_{i_1}\cdots q_{i_a}.
\endalign
$$
The coefficient $\la$ here is equal to
$\lan \Si^{2}_{(i_r)}  \vert 
\Si^\prime \vert 
\text{point} \ran_{q_{i_1}\cdots q_{i_a}},
$
where $\Si^\prime$ denotes the dual homology class to
$y_{i_1}^2\cdots \hat y_{i_r}^2 \cdots y_{i_a}^2$.
This means we are counting rational curves 
$f = (W_0\sub W_1\sub\cdots\sub W_n)$ with
$W_i(z)$ constant if $i\notin \{i_1,\dots,i_a\}$.  

Since no two of $i_1,\dots,i_a$ are
consecutive (mod $n$), $f$ may be
identified with a rational curve in a product 
$\C P^1\times\dots\times \C P^1$ of complex projective
lines.  It follows that $\la=0$, as required.
\qed
\enddemo

\subheading{\S 3 The periodic Toda lattice } 

We are now ready to prove that the relations in the algebra 
$QH^\sh \Fn$ are equal to the conserved quantities of the periodic Toda 
lattice: 

\proclaim{Theorem 3.1} For $0\le k\le n$, we have $QS^k_n = P^k_n$. 
\endproclaim

\no This will be an immediate consequence of Definition 2.4,
Corollary 2.5, and part (2) of the following Proposition 3.2. 
We use the notation
$O_n,P_n$ from the introduction.

\proclaim{Proposition 3.2}

\no(1) $O_n=
D_{n-1} D_{n-2} \cdots D_1 S_n$.

\no(2) $P_n=
D_n D_{n-1} \cdots D_1 S_n +
(-1)^{n+1} \frac{\q_1 \q_2 \cdots \q_n}{z} + z $. \endproclaim

\demo{Proof}
(1) Expanding $O_{k+1} = \det(L_{k+1} + \mu I)$ along the last row, we 
have
$$
\spreadlines{5\jot}
\align
O_{k+1}
&
=
(\x_{k+1}+\mu )
\eightpoint{
\left|
\matrix
\x_1+\mu & \cdots & \cdots & 0	\\
-1	& \cdots & \cdots & \cdots	\\
\vdots & \ddots & \ddots & \ddots	\\
0	& \cdots & \x_{k-1}+\mu & \q_{k-1}	\\
0	& \cdots & -1	& \x_k+\mu
\endmatrix
\right|
+
\left|
\matrix
\x_1+\mu & \cdots & \cdots & 0	\\
-1	& \cdots & \cdots & \cdots \\
\vdots & \ddots & \ddots	& \vdots \\
0	& \cdots & \x_{k-1}+\mu	& 0
	\\
0	& \cdots & -1	& \q_k
\endmatrix
\right|
}
\\
& =
(\x_{k+1}+ \mu) O_k +
\q_k (\frac{\partial}{\partial \x_k} O_k) \\ & =
(\x_{k+1}+ \mu) O_k +
\q_k\frac{\partial^2}{\partial \x_k \partial \x_{k+1}} \{ (\x_{k+1}+ \mu) 
O_k \} \\
& =
D_k \{ (\x_{k+1}+ \mu) O_k \} .
\endalign
$$
Since $(\x_j+\mu)D_i F = D_i\{ (\x_j+\mu)F \}$ for any polynomial $F$ if 
$j\ge i+2$,
part (1) now follows by induction.

(2) The only difference between $P_n$ and $O_n$ is that additional 
entries $-z$ and $\q_n/z$ appear in the top right and bottom left corners 
of the determinant. Expanding (partially) along the
last row, we see that $P_n$ is equal to 
$$
\spreadlines{5\jot}
\eightpoint{
\left|
\matrix
\x_1+\mu & \q_1	& \cdots
	& \cdots	& -z
	\\
-1	& \x_2+\mu & \cdots	& \cdots
	& 0	\\
\vdots & \ddots & \ddots	& \ddots
	& \vdots \\
0	& \cdots & \x_{n-2}+\mu & \q_{n-2}	& 0
	\\
0	& \cdots & -1	& \x_{n-1}+\mu 
& \q_{n-1}\\
0	& \cdots & 0	& -1
	& \x_n+\mu
\endmatrix
\right|
+
\left|
\matrix
\x_1+\mu & \q_1	& \cdots
	& \cdots	& -z
	\\
-1	& \x_2+\mu & \cdots	& \cdots
	& 0	\\
\vdots & \ddots & \ddots	& \ddots
	& \vdots \\
0	& \cdots & \x_{n-2}+\mu & \q_{n-2}	& 0
	\\
0	& \cdots & -1	& \x_{n-1}+\mu 
& \q_{n-1}\\
{\q_n}/{z} & \cdots & 0	& 0
	& 0
\endmatrix
\right|.
}
$$
Applying the same procedure to the right hand columns, the first 
determinant becomes
$$
\spreadlines{5\jot}
\eightpoint{
\left|
\matrix
\x_1+\mu & \q_1	& \cdots
	& \cdots	& 0 \\
-1	& \x_2+\mu & \cdots	& \cdots
	& 0	\\
\vdots & \ddots & \ddots	& \ddots
	& \vdots \\
0	& \cdots & \x_{n-2}+\mu & \q_{n-2}	& 0
	\\
0	& \cdots & -1	& \x_{n-1}+\mu 
& \q_{n-1} \\
0	& \cdots & 0	& -1
	& \x_n+\mu
\endmatrix
\right|
+
\left|
\matrix
\x_1+\mu & \q_1	& \cdots
	& \cdots	& -z
	\\
-1	& \x_2+\mu & \cdots	& \cdots
	& 0	\\
\vdots & \ddots & \ddots	& \ddots
	& \vdots \\
0	& \cdots & \x_{n-2}+\mu & \q_{n-2}	& 0
	\\
0	& \cdots & -1	& \x_{n-1}+\mu 
& 0 \\
0	& \cdots & 0	& -1
	& 0
\endmatrix
\right|
}
$$
while the second determinant (after expansion along the last row) 
becomes $$
\spreadlines{5\jot}
\eightpoint{
(-1)^{n \! + \! 1}\frac{\q_n}{z}
\left|
\matrix
\q_1	& \cdots & \cdots	& \cdots
	& 0	\\
\x_2 \! + \! \mu & \cdots & \cdots	& \cdots
	& 0	\\
\vdots & \ddots & \ddots	& \ddots
	& \vdots \\
0	& \cdots & \x_{n \! - \! 2} \! + \! \mu & \q_{n \! - \! 2}
	&
0	\\
0	& \cdots & -1 	& \x_{n \! - \! 1} 
\! + \! \mu & \q_{n \! - \! 1}
\endmatrix
\right|
+
(-1)^{n \! + \! 1}\frac{\q_n}{z}
\left|
\matrix
\q_1	& \cdots & \cdots	& \cdots
	& -z  \\
\x_2 \! + \! \mu & \cdots & \cdots	& \cdots
	& 0	\\
\vdots & \ddots & \ddots	& \ddots
	& \vdots \\
0	& \cdots & \x_{n \! - \! 2} \! + \! \mu & \q_{n \! - \! 2}
	&
0	\\
0	& \cdots & -1 	& \x_{n \! - \! 1} 
\! + \! \mu & 0
\endmatrix
\right|.
}
$$
The last term here is
$$
\spreadlines{5\jot}
(-1)^{n+1}\frac{\q_n}{z}(-1)^n(-z)
{\eightpoint
\left|
\matrix
\x_2+\mu & \cdots & \cdots	& \cdots
	\\
\ddots & \ddots & \ddots	& \vdots
	\\
\cdots & \cdots & \x_{n-2}+\mu & \q_{n-2}	\\
\cdots & \cdots & -1	& \x_{n-1}+\mu
\endmatrix
\right|
}
=
\q_n \frac{\partial^2}{\partial \x_n \partial \x_{1}}O_n. 
$$
Taking the sum of all four terms, we have
$$
\align
P_n&=
O_n +z
+ (-1)^{n+1}\frac{\q_n}{z} \q_1 \cdots \q_{n-1} + \q_n 
\frac{\partial^2}{\partial \x_n \partial \x_{1}}O_n \\ &=
D_n O_n + z + (-1)^{n+1}\frac{\q_1 \cdots \q_{n}}{z}. \endalign
$$
The required formula for $P_n$ follows from this and (1). 
\qed
\enddemo

\subheading{\S 4 Remarks}

Our computation of quantum products in \S 2 recovers
the result of Givental and Kim for $F_n$.  To see
this, we denote by  $I:F_n\to \Fn$ the inclusion of the fibre
$$
\pi_n^{-1}(H_+) =
\{\  \la H_+ \sub  W_1 \sub \cdots \sub W_{n-1} \sub H_+ \in \Fn\  \}
\cong F_n.
$$
Let $\hat y_i = I^\ast y_i$ and  $\hat x_i = I^\ast x_i$, where $1\le i\le n$ 
as usual.
Observe that $\hat y_0 = \hat y_n = 0$ now, and so
$\hat y_i = \hat x_1 + \dots + \hat x_i$ when $1\le i\le n$.  

In this situation we have evaluation maps
$$
\align
\hat{ev}_q&: \C[\y_1,\dots,\y_{n-1},\q_1,\dots,\q_{n-1}] 
\to QH^\ast F_n\\
\hat{ev}_c&:\C[\y_1,\dots,\y_{n-1},\q_1,\dots,\q_{n-1}] 
\to H^\ast F_n \otimes\C[q_1,\dots,q_{n-1}],
\endalign
$$
given by $\y_i\mapsto \hat y_i$ and $\q_i\mapsto q_i$, $1\le i\le n-1$.  
We denote
by $\hat V$ the $\C[\q_1,\dots,\q_{n-1}]$-submodule of 
$\C[\y_1,\dots,\y_{n-1},\q_1,\dots,\q_{n-1}]$ that is 
generated by elements of the form
$\x_{i_1} \cdots \x_{i_k}$, $1\le i_1 < \dots < i_k \le n$
where $\x_i=\y_i-\y_{i-1}$ (and $\y_0=\y_n=0$).
The analogue of Proposition 2.3 is then:

\proclaim{Proposition 4.1}  
On $\hat V$ we have

\no{(1)}
$\hat{ev}_q=
\hat{ev}_c  \delta_{n-1}\cdots \delta_1 $. 

\no{(2)}
$\hat{ev}_c=
\hat{ev}_q D_{n-1}\cdots D_1 $.
\endproclaim

\demo{Proof} The difference between the current situation and the 
situation
of Proposition 2.3 is that the expression $\hat y_n\circ \hat y_n$ 
contributes nothing to the quantum deformation.
Since the only other appearances of $\hat y_n$ are linear, the result
is the same as for  Proposition 2.3 but with the final $\de_n$ (and $D_n$) 
omitted.
\qed
\enddemo

Let 
$\hat{QS}_n = \sum_{i=0}^n \hat{QS}^i_n \mu^i = D_{n-1}\cdots D_1 S_n$.
By the argument of \cite{Si-Ti}, the coefficients $\hat{QS}^i_n$
are the defining relations for $QH^* F_n$ .

\proclaim{Corollary 4.2 (Givental and Kim)}  
For $0\le k\le n$, we have $\hat{QS}^k_n = O^k_n$.
\endproclaim

\demo{Proof}  Proposition 4.1 and part (1) of Proposition 3.2.
\qed
\enddemo

A comment is necessary on our \ll differential operator formulae\rr for the
quantum relations in the case of $F_n$ (Proposition 4.1). When we
discovered these formulae we believed (naively) that they were new. However,
after completing our calculations, we became aware of (i) the papers 
\cite{Sa-Ko}, \cite{Wo} in which similar formulae were given for the
conserved quantities of the Toda lattice, and (ii) the paper \cite{Ci2}
(containing full details of the results  announced in \cite{Ci1}) in which similar
formulae were obtained for the quantum relations of $F_n$ as a consequence
of a general theory of quantum Schubert calculus.

Returning to the infinite-dimensional case,
it should be said that the relation between $QH^\sh \Fn$ 
and the periodic Toda lattice is a plausible extension of the formula
of Givental and Kim, in view of the
following two facts:

\no(a) Formally, the open Toda lattice may be obtained from the 
periodic Toda lattice by setting $\q_n=0$.

\no(b) For a finite-dimensional fibre bundle, formula (2.17) of \cite{As-
Sa}
says that the quantum cohomology of the fibre should be obtained by 
dividing 
the vertical quantum cohomology of the  total space by the cohomology 
(in positive dimensions) of the base.   Applying this
to the bundle $\pi_n: \Fn \to \Gn$ amounts to setting $q_n=0$
(to obtain the vertical quantum cohomology) and then $y_n=0$.
Our formula is consistent with this procedure.

In this paper we have focused attention on $QH^\sh \Fn$, \ll the
subalgebra of $QH^\ast \Fn$ generated by two-dimensional classes\rrr,
in accordance with the philosophy of \cite{Gi-Ki}, \cite{Au}.  However,
as we remarked in the introduction, it is not {\it a priori} clear whether
$QH^\sh \Fn$ is the same as $H^\sh\Fn \otimes \C[q_1,\dots,q_{n}]$,
i.e. whether the latter space is closed under quantum multiplication. 
The smaller subspace $H^\ast F_n \otimes \C[q_1,\dots,q_{n}]$ is
in fact more appropriate from the point of view of the periodic Toda lattice,
and its use will render assumption (2) of the introduction unnecessary.
(Note that $H^\ast F_n$ is a subalgebra of the ordinary cohomology algebra
$H^\ast \Fn$ because $\Fn$ is diffeomorphic to $F_n \times \Gn$.) It
may be shown by the methods of this paper that 
$H^\ast F_n \otimes \C[q_1,\dots,q_{n}]$ is closed under our
hypothetical quantum product, and hence that we have an isomorphism
of algebras
$$
H^\ast F_n \otimes \C[q_1,\dots,q_{n}]\cong
\C[\x_1,\dots,\x_n,\q_1,\dots,\q_n]/
\lan QS^0_n, QS^1_n,\dots,QS^{n-1}_n \ran.
$$
Furthermore, this algebra is the \ll coordinate ring\rr of
a spectral cover (in the sense of \cite{Au}), which is in
turn exactly the zero level set of the conserved quantities
of the periodic Toda lattice.

Finally, we point out that the quantum cohomology calculation of
Lemma 2.6 --- the main ingredient of Propositions
2.3 and 4.1 --- amounts to an inductive procedure whereby certain 
quantum
products in a full flag manifold are reduced to quantum products in  
products of flag manifolds of lower rank.  
This seems likely to work quite generally for
flag manifolds of the form $G/B$ or $LG/B$.

\eightpoint \it

\no MG:

\no  Department of Mathematics
\newline
Graduate School of Science
\newline
Tokyo Metropolitan University
\newline
Minami-Ohsawa 1-1, Hachioji-shi
\newline
Tokyo 192-0397, Japan

\no martin\@comp.metro-u.ac.jp

\no TO:

\no  Department of Mathematics
\newline
Graduate School of Science and Engineering
\newline
Tokyo Institute of Technology
\newline
Okayama 2-12-1, Meguro-ku
\newline
Tokyo 152-8551, Japan

\no otofuji\@math.titech.ac.jp


\newpage

\Refs   
     
\widestnumber\key{Fo-Ge-Po}

\ref
\key  As-Sa
\by A. Astashkevich and V. Sadov
\paper Quantum cohomology of partial flag manifolds
\jour Comm. Math. Phys.
\yr 1995
\vol 170
\pages 503--528
\endref

\ref
\key At
\by M.F. Atiyah
\pages 437--451
\paper Instantons in two and four dimensions
\yr 1984
\vol 93
\jour Comm. Math. Phys.
\endref

\ref 
\key  Au
\by M. Audin
\paper Symplectic geometry in Frobenius manifolds
and quantum cohomology
\jour J. Geom. and Physics
\yr 1998
\vol 25
\pages 183--204
\endref

\ref
\key  Be
\by A. Bertram
\paper Quantum Schubert calculus
\jour Advances in Math.
\yr 1997
\vol 128
\pages 289--305
\endref

\ref\key  Bo \by R. Bott
\paper  The space of loops on a Lie group
\yr 1958
\vol 5
\jour Michigan Math. J.
\pages 35--61
\endref

\ref
\key Ca
\by J.B. Carrell
\paper Vector fields, flag varieties, and Schubert calculus
\inbook Proceedings of the Hyderabad Conference on 
Algebraic Groups
\ed S. Ramanan
\publ Manoj Prakashan, Madras
\yr 1991
\pages 23--57
\endref

\ref 
\key  Ci1
\by I. Ciocan-Fontanine
\paper Quantum cohomology of flag varieties
\jour Internat. Math. Res. Notices
\yr 1995
\vol 6
\pages 263--277
\endref

\ref 
\key  Ci2
\by I. Ciocan-Fontanine
\paper The quantum cohomology ring of flag varieties
\jour Trans. Amer. Math. Soc.
\yr 1999
\vol 351
\pages 2695--2729
\endref

\ref
\key  Fo-Ge-Po
\by S. Fomin, S. Gelfand, and A. Postnikov
\paper Quantum Schubert polynomials
\jour J. Amer. Math. Soc.
\yr  1997
\vol 10
\pages 565--596
\endref

\ref
\key  Fr
\by D.S. Freed
\paper The geometry of loop groups
\jour J. Differential Geom.
\yr 1988
\vol 28
\pages 223--276
\endref

\ref
\key  Gi1
\by A.B. Givental
\paper Equivariant Gromov-Witten invariants
\jour Internat. Math. Res. Notices
\yr 1996
\vol 13
\pages 613--663
\endref

\ref
\key  Gi2
\by A. Givental
\paper Stationary phase integrals, quantum Toda lattices, flag
manifolds and the mirror conjecture
\inbook Topics in Singularity Theory 
\bookinfo AMS Translations 180
\eds A. Khovanskii, A. Varchenko, and V. Vassiliev
\yr 1997
\publ Amer. Math. Soc.
\pages  103--115
\endref

\ref 
\key  Gi-Ki
\by  A. Givental and B. Kim
\paper Quantum cohomology of flag manifolds and Toda lattices
\jour Comm. Math. Phys.
\yr 1995
\vol 168
\pages 609--641
\endref

\ref
\key  Gu
\by E. Gutkin
\paper Schubert calculus on flag varieties of Kac-Moody groups
\jour Algebras, Groups and Geometries
\yr 1986
\vol 3
\pages 27-59
\endref

\ref
\key Ha
\by Z. Haddad
\paper A Coxeter group approach to Schubert varieties
\inbook Infinite-dimensional groups with applications
\bookinfo Math. Sci. Res. Inst. Publ. 4
\publ Springer
\ed V. Kac
\yr 1985
\pages 157--165
\endref

\ref
\key  Ki1
\by B. Kim
\paper On equivariant quantum cohomology
\jour Internat. Math. Res. Notices
\yr  1996
\vol 17
\pages 841--851
\endref

\ref
\key  Ki2
\by B. Kim
\paper Quantum cohomology of flag manifolds $G/B$ and
quantum Toda lattices
\jour Ann. of Math.
\yr 1999
\vol 149
\pages 129--148
\endref

\ref
\key  Ko1
\by  B. Kostant
\paper Flag manifold quantum cohomology, the Toda lattice, 
and the representation with highest weight $\rho$
\jour Selecta Math. (N.S.)
\yr 1996
\vol 2
\pages 43--91
\endref

\ref
\key  Ko2
\by B. Kostant
\paper Quantum cohomology 
of the flag manifold as an algebra of rational functions on a
unipotent algebraic group
\inbook   Deformation theory and symplectic 
geometry (Ascona, 1996)
\eds  D. Sternheimer, J. Rawnsley and S. Gutt
\bookinfo Mathematical Physics Studies, 20
\publ Kluwer
\yr  1997
\pages  157--175
\endref

\ref
\key  Ko-Ku
\by B. Kostant and S. Kumar
\paper The nil Hecke ring and cohomology of $G/P$
for a Kac-Moody group $G$
\jour Advances in Math.
\yr 1986
\vol 62
\pages 187--237
\endref

\ref
\key  Lu
\by P. Lu
\paper A rigorous definition of fiberwise quantum cohomology
and equivariant quantum cohomology
\jour Comm. Anal. Geom.
\yr 1998
\vol 6
\pages 511--588
\endref

\ref
\key Ol-Pe
\by M.A. Olshanetsky and A.M. Perelomov
\paper  Integrable systems and finite-dimensional Lie algebras
\inbook  Dynamical Systems VII
\bookinfo Encylopaedia of Mathematical Sciences 16
\eds  V.I. Arnol'd and S.P. Novikov
\publ Springer
\yr 1994
\pages 87--116
\endref

\ref
\key Pe
\by A.M. Perelomov
\book Integrable Systems of Classical Mechanics and Lie Algebras
\publ Birkh\"auser
\yr 1990
\endref

\ref
\key  Pr-Se
\by A.N. Pressley and G.B. Segal
\book  Loop Groups
\publ Oxford Univ. Press 
\yr1986
\endref

\ref
\key Re-Se
\by A.G. Reyman and M.A. Semenov-Tian-Shansky
\paper  Group-theoretical methods in the theory of
finite-dimensional integrable systems
\inbook  Dynamical Systems VII
\bookinfo Encylopaedia of Mathematical Sciences 16
\eds  V.I. Arnol'd and S.P. Novikov
\publ Springer
\yr 1994
\pages 116--225
\endref

\ref
\key Sa-Ko
\by K. Sawada and T. Kotera
\paper Integrability and a solution for the one-dimensional $N$-particle system
with inversely quadratic pair potentials
\jour J. Phys. Soc. Japan 
\vol 39 
\yr 1975
\pages 1614--1618
\endref

\ref
\key  Si-Ti
\by B. Siebert and G. Tian
\paper On quantum cohomology rings of Fano manifolds and
a formula of Vafa and Intriligator
\paperinfo alg-geom/9403010
\jour Asian J. Math.
\yr 1997
\vol 1
\pages  679--695
\endref

\ref
\key Wo
\by S. Wojciechowski
\paper Involutive set of integrals for completely integrable
many-body problems with pair interaction
\jour  Lett. Nuovo Cimento
\vol 18
\yr 1977
\pages 103--107
\endref

\endRefs
\enddocument